\documentclass{amsart}
\usepackage{amsfonts}
\usepackage{amsmath,amsthm}
\usepackage{amsfonts,amssymb}
\usepackage{mathrsfs}
\usepackage{enumerate}
\usepackage{bm}

\newtheorem{thm}{Theorem}[section]

\theoremstyle{remark}
\newtheorem{rem}[thm]{Remark}

\numberwithin{equation}{section}

\def\az{\alpha}
\def\bz{\beta}
\def\gz{\gamma}

\def\dz{\delta}

\def\vz{\varepsilon}
\def\lz{\lambda}

\def\sz{\sigma}

\def\vr{\varphi}
\def\ka{\mathcal{K}}
\def\ze{\zeta}

\def\s{{\bf s}}

\def\va{\varepsilon}

\def\t{{\bf t}}
\def\s{{\bf s}}
\def\x{{\bf x}}

\def\y{{\bf y}}

\def\B{\Bbb B}
\def\N{\mathbb N}

\begin{document}

\title[] {Average case tractability of additive random fields with Korobov kernels}

\author[]{Jia Chen, Heping Wang} \address{School of Mathematical Sciences, Capital Normal
University, Beijing 100048,
 China.}
\email{ jiachencd@163.com;\ \  \ wanghp@cnu.edu.cn.}


\begin{abstract} We  investigate  average case tractability of
approximation of
 additive random fields  with marginal random processes
corresponding to the Korobov kernels  for  the non-homogeneous
case.  We use  the absolute error criterion (ABS)   or  the
normalized error  criterion (NOR). We  show that  the problem   is
always polynomially tractable for ABS or NOR, and  give sufficient
and necessary conditions for strong polynomial tractability  for
ABS or NOR.

\end{abstract}

\maketitle
\input amssym.def

\section{Introduction}

Let $X_j(t),\, t\in [0,1],\,j\in \N$ be a sequence of independent
random processes, where  $\N=\{1,2,\cdots\}$.
 Suppose that every  random element $X_j(t)$
has zero mean and a covariance function $\ka^{X_j}(t,s),\,t,\, s\in [0,1]$. Let $K^{X_j}$ be the covariance operator with the
covariance kernel $\ka^{X_j}$ for every $j\in \N$. Then for any $f\in L_2([0,1])$ and $t\in [0,1]$, we have
$$ K^{X_j}f(t)=\int_{[0,1]}\ka^{X_j}(t,s)f(s)ds.$$
 We consider the random field
\begin{equation}\label{1.1}Y_d (\t) :=\sum_{j=1}^dX_j(t_j),\, \ \ \t=(t_1 ,\cdots,t_d )\in [0,1]^d,\end{equation}
with the following zero-mean covariance function
\begin{equation}\label{1.2}\ka^{Y_d}(\t,\s)=\sum_{j=1}^d\ka^{X_j}(t_j,s_j ),\end{equation}
and covariance operator
\begin{equation}\label{1.3}K^{Y_d}f(\t)=\int_{[0,1]^d}\ka^{Y_d}(\t,\s)f(\s)d\s,\quad  \end{equation} where $f\in
L_2([0,1]^d), \ \t=(t_1 ,\cdots, t_d ),\, \s=(s_1 ,\cdots, s_d ),\
\t,\s\in [0,1]^d.$ This random field is called an additive random
field. There are many papers which investigated this random field,
see \cite{CL, KZ1, KZ2, LZ}.

In this paper, we investigate the average case approximation of
$Y_d(\t),\, \t\in[0,1]^d$, as a random element of the space
$L_2([0,1]^d)$ equipped with  inner product $(\cdot,\cdot)_{2,d}$
and norm $\|\cdot\|_{2,d}$, by a finite rank random field.

The $n$th minimal 2-average case error, for $n\ge 1$, is defined by
$$e^{Y_d}(n):=\inf\big\{\big(\mathbb{E}\|Y_d-Y_d^{(n)}\|^2_{2,d}\big)^{1/2}\,:\,Y_d^{(n)}\in \mathcal{A}^{Y_d}_n\big\}.$$
Here $\mathcal{A}^{Y_d}_n$ is the class of all linear algorithms
with rank $\le n$ defined by
$$\mathcal{A}^{Y_d}_n:=\big\{\sum_{m=1}^n(Y_d ,\vr_m )_{2,d}\,\psi_m : \vr_m,\psi_m\in L_2( [ 0,1
]^d)\big\}.$$

The so-called initial average case error is given by
$$e^{Y_d}(0):=\big(\mathbb{E}\|Y_d\|^2_{2,d}\big)^{1/2}.$$

We use either  the absolute error criterion (ABS) or  the
normalized error criterion (NOR). For any $\va\in (0,1)$,  we
study information complexity $n^{Y_d,Z} (\va)$ of approximation of
the random fields $Y_d,\,d\in \N$ defined by
$$n^{Y_d,Z} (\va) := \min \{ n\in \N : e^{Y_d}(n)\le \va {\rm CRI}_d \} , $$
where  $Z\in\{{\rm ABS},{\rm NOR}\}$,
\begin{equation*}
{\rm CRI}_d:=\left\{\begin{matrix}
 & \; 1, \; \qquad\text{ for} \ Z={\rm ABS}, \\
 & e^{Y_d}(0),\quad \text{for} \ Z={\rm NOR}.
\end{matrix}\right.
\end{equation*}

Let $Y=\{Y_d\}_{d\in\Bbb N}$. First we  consider average case
tractability of $Y$. Various notions of tractability have been
discussed for multivariate problems. We recall
 some of the basic tractability notions (see \cite{NW1, NW2, NW3, S1}).

For $Z\in \{{\rm ABS,\, NOR}\}$,
say $Y$ is

$\bullet$   {\it strongly polynomially tractable (SPT)}   iff
there exist non-negative numbers $C$ and $p$ such that for all
$d\in \Bbb N,\ \va \in (0,1)$,
\begin{equation*}
n^{Y_d,Z} (\va) \leq C(\va ^{-1})^p;
\end{equation*}
The infimum of $p$ satisfying the above inequality is called the exponent
of strong polynomial tractability and is denoted by $p^{\rm str-avg}$.

 $\bullet$  {\it polynomially tractable
(PT)} iff there exist non-negative numbers $C,\, p$ and $q$ such
that for all $d\in \Bbb N, \ \va \in(0,1)$,
\begin{equation*}
n^{Y_d,Z} (\va) \leq Cd^q(\va ^{-1})^p.
\end{equation*}

$\bullet$   {\it quasi-polynomially tractable (QPT)} iff there
exist two constants $C,t>0$ such that for all $d\in \Bbb N, \ \va
\in(0,1)$,
\begin{equation*}
n^{Y_d,Z} (\va) \leq C\exp(t(1+\ln\va ^{-1})(1+\ln d));
\end{equation*}

$\bullet$ {\it uniformly weakly tractable (UWT)} iff for all $s, t>0$,
\begin{equation*}
\lim_{\varepsilon ^{-1}+d\rightarrow \infty }\frac{\ln n^{Y_d,Z} (\va) }{(\va ^{-1})^{s }+d^{t }}=0;
\end{equation*}

$\bullet$  {\it weakly tractable (WT)} iff
\begin{equation*}
\lim_{\va ^{-1}+d\rightarrow \infty }\frac{\ln n^{Y_d,Z} (\va) }{\va
^{-1}+d}=0.
\end{equation*}

This paper is devoted to studying average case  tractability of
the additive random field $Y=\{Y_d\}_{d\in \N}$ under ABS and NOR.
For additive random fields similar problems were investigated in
\cite{HWW, LZ, LZ2, WW} in various settings for the homogeneous
case and in \cite{KZ1, KZ2} for  the non-homogeneous case. Here,
the homogeneous case means that  approximated additive random
fields are constructed (in a special way) from copies of one
marginal process, while   the non-homogeneous case  means that the
random fields are composed of a whole sequence of marginal random
processes with generally different covariance functions.
Specifically, the authors in \cite{KZ1} obtained the growth of
$n^{Y_d,{\rm NOR}}(\va)$ for arbitrary fixed $\vz\in(0,1)$ and
$d\to\infty$ for the non-homogeneous case and gave application to
the additive random fields with marginal random processes
corresponding to the Korobov kernels.

 It
should be noted, however, that all these works deal only with NOR.
In this paper, we consider average case  tractability of the
problem $\B=\{\B_d\}_{d\in\N}$ of the additive random fields  with
marginal random processes corresponding to the Korobov kernels
under ABS and NOR. We shall show that  the problem $\B$ is always
polynomially tractable for ABS or NOR. Obviously, PT implies all
QPT, UWT, WT. We also give sufficient and necessary conditions for
which $\B$ is SPT for ABS or NOR.

The paper is organized as follows. In Section 2 we give
preliminaries about the additive random fields  with marginal
random processes corresponding to the Korobov kernels and
introduce main results, i.e., Theorems 2.1-2.3.  Section 3 is
devoted to proving Theorems 2.1-2.3.

\section{Preliminaries and main results}

Let us consider  the sequence of additive random fields $Y_d(\t)$,
$\t \in [0,1]^d$, $d\in\N$, defined by \eqref{1.1}. Let
$\big\{(\lz _k^{X_j},\,\psi _k^{X_j})\big\}_{k\in\N}$ be the
sequence of eigenpairs of the covariance operator $K^{X_j}$ of
$X_j$. Under the additive structure \eqref{1.2}, the eigenvalues
$\lz_{d,k},\,k\in\N$, are generally unknown or not easily depend
on $\lz _k^{X_j},\,k\in\N,\,j=1,2,\cdots,d$. However, under the
following condition, we can explicitly describe the eigenvalues
$\lz_{d,k},\,k\in\N$.

For every $j\in\N$ there exist
$\psi_0\in\{\psi_1^{X_j},\psi_2^{X_j},\cdots,\psi_k^{X_j},\cdots\}$
such that $\psi_0(t)=1$ for all $t\in[0,1]$. We denote by
$\bar{\lz}_0^{X_j}$ the eigenvalues corresponding to eigenvector
1. Let $\{\bar{\lz}_k^{X_j}\}_{k\in\N}$ and
$\{\bar\psi_k^{X_j}\}_{k\in\N}$ be the non-increasing sequence of
the remaining eigenvalues and the corresponding sequence of
eigenvectors of $K^{X_j}$, respectively. It is known that the
family of eigenvectors $${1}\cup\{\bar
{\psi}_k^{X_j}(t_j):k\in\N,\,j=1,\cdots,d\},\quad \t=(t_1,\cdots
,t_d)\in[0,1]^d,$$ is an orthogonal system in $L_2([0,1]^d)$ for
every $d\in\N$, see \cite{KZ1}. Hence the identical 1 is an
eigenvector of $K^{Y_d}$ with the eigenvalue
$\bar{\lz}_0^{Y_d}:=\sum_{j=1}^d\bar{\lz}_0^{X_j}$, and the pairs
$(\bar{\lz} _k^{X_j},\,\bar{\psi} _k^{X_j})$, for all $k\in\N$ and
$j=1,\cdots,d$, are the remaining eigenpairs of $K^{Y_d}$.

 Let
$\{\lz _{d,j}\}_{j\in \Bbb N}$ and $\{\psi _{d,j}\}_{j\in \Bbb N}$
be the non-increasing sequence of the eigenvalues and the
corresponding sequence of eigenvectors of  $K^{Y_d}$ defined by
\eqref{1.3}. Then the average case information complexity
$n^{Y_d,Z} (\va)$ can also be described in terms of eigenvalues
$\lz_{d,j},\,j\in\N$ of $K^{Y_d}$ by
\begin{equation}\label{1.4}n^{Y_d,Z} (\va) := \min \{ n\in \N : \sum_{j=n+1}^\infty\lz_{d,j}\le \va^2 {\rm CRI}_d^2 \} , \end{equation}
where
\begin{equation*}
{\rm CRI}_d=\left\{\begin{matrix}
 & \quad 1, \,\, \qquad\quad\text{ for} \ Z={\rm ABS}, \\
 & \big(\sum_{j=1}^\infty\lz_{d,j}\big)^{1/2},\quad \text{for} \ Z={\rm NOR},
\end{matrix}\right.
\end{equation*}
see \cite{NW1}.

Particularly, we study  additive random fields with marginal random processes  corresponding to the Korobov kernels.
Let $B_{\az,\bz,\sz}(x)$, $x\in[0,1]$ be a zero-mean random field with the following covariance function
\begin{equation*}\ka_{\az,\beta,\sz}(x,y)=\az+2\bz\sum_{k=1}^\infty k^{-\sz}\cos(2\pi k(x-y)),\quad x,\, y\in[0,1].\end{equation*}
Here $\az\ge 0$, $\bz>0$ and $\sz>1$. Let $K^{B_{\az,\bz,\sz}}$ be
the covariance operator with kernel
 $\ka_{\az,\bz,\sz}$ of $B_{\az,\bz,\sz}$, and for any
 $L_2([0,1])$ and $x\in[0,1]$, $$K^{B_{\az,\bz,\sz}}f(x):=\int_{[0,1]}\ka_{\az,\bz,\sz}(x,y)f(y)dy.$$

The eigenpairs of the covariance operator $K^{B_{\az,\bz,\sz}}$
are known, see \cite{NW1}. The identical 1 is an eigenvector of
$K^{B_{\az,\bz,\sz}}$ with the eigenvalue
$\bar{\lz}_0^{B_{\az,\bz,\sz}}=\az$. The other eigenpairs of
$K^{B_{\az,\bz,\sz}}$ have the following form:
$$\bar{\lz}_{2k-1}^{B_{\az,\bz,\sz}}=\bar{\lz}_{2k}^{B_{\az,\bz,\sz}}=\frac{\bz}{k^{\sz}},\, \bar{\psi}_{2k-1}^{B_{\az,\bz,\sz}}(x)=e^{-i2\pi kx},\, \bar{\psi}_{2k}^{B_{\az,\bz,\sz}}(x)=e^{i2\pi kx},$$
for any $k\in \N,\, x\in[0,1]$.

Suppose that $B_j(x),\, x\in [0,1], \, j\in\Bbb N$ is a sequence
of independent zero-mean random fields with covariance functions
$\ka_{\az_j,\bz_j,\sz_j}$, respectively. Let
$\B_d(\x)=\sum_{j=1}^dB_j(x_j)$, $\x\in [0,1]^d,\, d\in \N$, be
the sequence of zero-mean random fields with the covariance
functions
\begin{equation*}\ka^{\B_d}(\x,\y)=\sum_{j=1}^d\ka_{\az_j,\bz_j,\sz_j}(x_j,y_j),\end{equation*}
where $\x=(x_1,x_2,\cdots,x_d)\in[0,1]^d$,
$\y=(y_1,y_2,\cdots,y_d)\in[0,1]^d$,  and the parameters $\az_j\ge
0$, $\bz_j>0$ for all $j\in\N$, and $\inf_{j\in \N}\sz_j>1$.

Let $K^{\B_d}$ be the covariance operator of $\B_d$.  We have
\begin{equation*}K^{\B_d}f(\x)=\int_{[0,1]^d}\sum_{j=1}^d\ka_{\az_j,\bz_j,\sz_j}(x_j,y_j)f(y_j)d\y,\, \end{equation*}
for any $f(\x)\in L_2([0,1]^d)$ and $\x \in[0,1]^d$.
Then the identical 1 is an eigenvector of $K^{\B_d}$ with
the eigenvalue
\begin{equation}\label{2.1}\bar{\lz}^{\B_d}_0=\sum_{j=1}^d\az_j\end{equation} and the remaining eigenvalues and eigenvectors are
\begin{equation}\label{2.2}\bar{\lz}^{\B_d}_{2k-1}=\bar{\lz}^{\B_d}_{2k}=\frac{\bz_j}{k^{\sz_j}},\ \ k\in \Bbb N,\ j=1,\dots,d,\end{equation}
 and $$\bar {\psi}_{2k-1}^{\B_d}(x_j)=e^{-i2\pi kx_j},\quad \bar {\psi}_{2k}^{\B_d}(x_j)=e^{i2\pi kx_j},\  \ \x\in[0,1]^d, \ k\in \Bbb N,\ j=1,\dots,d, $$
 respectively.  Let $\{\lz_{d,j}\}_{j=1}^\infty$
 be the sequence of non-increasing rearrangement of the eigenvalues of $K^{\B_d}$. Then we have
\begin{equation}\label{2.3}\sum_{j=1}^\infty \lz_{d,j}=\sum_{j=1}^d\az_j+2\sum_{j=1}^d\sum_{k=1}^\infty\frac{\bz_j}{k^{\sz_j}}
=\sum_{j=1}^d\az_j+2\sum_{j=1}^d\bz_j\ze(\sz_j),\end{equation} and
for any $\tau\in(0,1)$ and $\tau \sz_j>1, \ j=1,\dots,d$,
\begin{equation}\label{2.4}
\sum_{j=1}^\infty \lz_{d,j}^\tau=\big(\sum_{j=1}^d\az_j\big)^\tau+2\sum_{j=1}^d\sum_{k=1}^\infty\frac{\bz_j^\tau}{k^{\tau\sz_j}}=
\big(\sum_{j=1}^d\az_j\big)^\tau+2\sum_{j=1}^d\bz_j^\tau\ze(\tau\sz_j),
\end{equation}
where $\ze(p)=\sum_{k=1}^\infty k^{-p}, \, p>1$, is the Riemann zeta-function.

In the sequel  we always assume that the sequences
$\az_j,
 \,\beta_j,\,\sz_j$ satisfy
\begin{equation}\label{2.002} \az_j\ge0,\ \ \ \
1\ge\bz_1\ge\bz_2\ge\cdots>0, \ \
 \quad 1<\sz_1\le\sz_2\le\cdots.
\end{equation}

In this paper, we consider the tractability of the problem
$$\B=\{\B_d\,\,: \,\,L_2([0,1]^d)\to L_2([0,1]^d)\}$$ under ABS and NOR, where the sequences $\az_j,
 \,\beta_j,\,\sz_j$ satisfy \eqref{2.002}. Our main results can be formulated as follows.

\begin{thm}\label{th2.1}
Let the sequences $\az_j,\, \bz_j,\, \sz_j$ satisfy \eqref{2.002}.
Then the problem $\B$

(i) is always PT for ABS or NOR;

(ii) is SPT   for ABS iff
\begin{equation}\label{2.4-1}
A_*:=\underset{d\to\infty}{\underline{\lim}}\frac{\ln
\frac1{\bz_d}}{\ln d}>1.
\end{equation}
The exponent of SPT  is
\begin{equation}\label{2.4-2}
p^{{\rm str-avg}}=\max\big\{\frac2{A_*-1},\frac2{\sz_1-1}\big\}.
\end{equation}
\end{thm}

\

In order to investigate the strong polynomial  tractability of the
problem $\B$ for NOR, we consider two cases.

\begin{thm}\label{th2.2}
Let $\az_j,\, \bz_j,\, \sz_j$ satisfy \eqref{2.002} and
$$0\le\az_j\le c\bz_j,\quad \text{for all}\quad j\in \N,$$ where $c>0$ is a constant. Then for NOR the problem $\B$ is SPT iff
\begin{equation*}
A_*:=\underset{d\to\infty}{\underline{\lim}}\frac{\ln
\frac1{\bz_d}}{\ln d}>1.
\end{equation*}
The exponent of SPT  is
\begin{equation*}
p^{{\rm str-avg}}=\max\big\{\frac2{A_*-1},\frac2{\sz_1-1}\big\}.
\end{equation*}
\end{thm}

\begin{thm}\label{th2.3}
Let  $\az_j,\, \bz_j,\, \sz_j$ satisfy \eqref{2.002}. Further
assume that
 $$0<r_1\le \dots\le r_j\le r_{j+1}\le \cdots,\ \  {\rm where}\   r_j=
 \frac{\az_j}{\beta_j},\ j\in \Bbb N,$$ and there exists a constant $c>0$ such that for every $d\in\Bbb N$,  $${\sum_{j=1}^d\az_j}\le c {d\az_d}.$$
 Then   the problem $\B$ is SPT for NOR iff\begin{equation}\label{2.5}
B_*:=\underset{d\to\infty}{\underline{\lim}}\frac{\ln
\frac{\az_d}{\bz_d}}{\ln d}>0.
\end{equation} The exponent of SPT  is
\begin{equation}\label{2.6}
p^{{\rm str-avg}}=\max\big\{\frac2{B_*},\frac2{\sz_1-1}\big\}.
\end{equation}

In particular,  if $\az_j=1$ for $j\in\Bbb N$, then   the problem
$\B$ is SPT for NOR iff
\begin{equation*}
A_*:=\underset{d\to\infty}{\underline{\lim}}\frac{\ln
\frac1{\bz_d}}{\ln d}>0.
\end{equation*}

\end{thm}

\begin{rem}In \cite{KZ1}, Khartov and Zani investigated
$n^{\B_d,{\rm NOR}}(\va)$ for arbitrary fixed $\vz\in(0,1)$ and
$d\to\infty$  of the above problem $\B$ with the
 parameters \begin{equation}\label{2.5-1}\bz_j\sim
cj^{-s},\quad \az_j/\bz_j\to r,\quad \sz_j\to +\infty,\quad
j\to\infty,\end{equation}  where $c>0$, $s>0$, and $0\le
r\le\infty$. They  obtained the following results.

(1) For $c>0$ and either $s>1$, $0\le r\le \infty$ or $0\le
s\le1$, $r=\infty$. Then
$$\sup_{d\in\N}n^{\B_d,\text{NOR}}(\va)<\infty,\quad\text{for every}\quad \vz\in(0,1).$$

(2)  For $c>0$, $0\le s< 1$ and $0\le r<\infty$,
$$n^{\B_d,\text{NOR}}(\va)\sim2Q(\va)\cdot d,\quad d\to \infty,\quad\text{for every }\vz
\in(0,\vz_0),$$
where
$$Q(\va)=\big(1-(\va/\va_0)^2\big)^{\frac1{1-s}},\ \ \va_0=(1+r/2)^{-1/2}. $$

(3)  For $c>0$, $s=1$ and $0\le r< \infty$, then
$$\ln n^{\B_d,\text{NOR}}(\va)=\big(1-(\va/\va_0)^2\big)\cdot\ln d +o(\ln d),\quad d\to\infty, \quad \vz\in(0,\va_0).$$

We remark that in some sense the above results give the
tractability results of  the problem $\B$ for NOR under the
condition \eqref{2.5-1}. Specifically, the above result (1)
corresponds to SPT of $\B$ for NOR, and results (2) and (3) relate
to PT of $\B$ for NOR, but more explicitly.  Comparing with
\cite{KZ1}, we obtain the tractability of the problem $\B$ with
general parameters $\az_j, \ \beta_j,\ \sz_j$ satisfying
\eqref{2.002} for ABS and NOR. Also  we get the exponent of SPT
for ABS and NOR.

\end{rem}

\begin{rem}\label{rem2.1} For the above problem $\B=\{\B_d\}_{d\in\Bbb N}$, let
$\widetilde\B_d$ be the additive random fields with marginal random
processes corresponding to the Korobov  covariance functions
$$\ka^{\widetilde\B_d}(\x,\y)=\sum_{j=1}^d\ka_{\widetilde\az_j,\widetilde\bz_j,\widetilde\sz_j}(x_j,y_j),\ \ \x,\y\in[0,1]^d,$$
with parameters $\widetilde \az_j,\widetilde\beta_j,\widetilde\sz_j$
satisfying
$$\widetilde \az_j=0,\ \ \widetilde\beta_j=\beta_j,\ \ \widetilde\sz_j=\sz_j,\
\ j\in\Bbb N.$$

 Let $\{{\lz}_{d,j}\}_{j=1}^\infty$ and $\{\widetilde{\lz}_{d,j}\}_{j=1}^\infty$  be the sequences of non-increasing rearrangement
  of the eigenvalues of the covariance operators $K^{{\B}_d}$ and
  $K^{\widetilde{\B}_d}$, respectively. Then  for some $j_0\in\Bbb N$,
 $$\lz_{d,j_0}=\sum_{j=1}^d\az_j, $$ and $$\widetilde
 \lz_{d,j}=\lz_{d,j}\ \  {\rm if  }\ \ j<j_0,\ \ \ {\rm and} \ \ \ \widetilde
 \lz_{d,j}=\lz_{d,j+1} \ \ {\rm if \ \  }j\ge j_0.$$
From the above equalities we obtain that $$\ \ \ \lz_{d,j}\ge
 \widetilde \lz_{d,j}\ge \lz_{d,j+1}, \ \ \ j\in\Bbb N.$$
 It follows from
\eqref{1.4} and the inequality $\lz_{d,j}\ge
 \widetilde \lz_{d,j}$ that for $\va\in(0,1)$,
 $$n^{\B_d,\text{ABS}}(\va)\ge n^{\widetilde{\B}_d,\text{ABS}}(\va).$$

On the other hand, we note that for $\va\in(0,1)$,
$$\sum_{j=n_0+1}^\infty\widetilde{\lz}_{d,j}\le \vz^2, $$where
$n_0=n^{\widetilde{\B}_d,\text{ABS}}(\va)$. This gives that
$$\sum_{j=n_0+2}^\infty{\lz}_{d,j}\le \sum_{j=n_0+1}^\infty\widetilde{\lz}_{d,j}\le \vz^2,
$$which means that $$n^{\B_d,\text{ABS}}(\va)\le
n_0+1=n^{\widetilde{\B}_d,\text{ABS}}(\va)+1.$$

 Hence we have for $\vz\in(0,1)$,
 $$n^{\widetilde{\B}_d,\text{ABS}}(\va)\le n^{\B_d,\text{ABS}}(\va)\le n^{\widetilde{\B}_d,\text{ABS}}(\va)+1.$$
It follows that  the problems $$\B_d: L_2([0,1]^d)\to
L_2([0,1]^d)\quad\text{and}\quad \widetilde{\B}_d: L_2([0,1]^d)\to
L_2([0,1]^d),\quad d\in\N,$$ have the same tractability for ABS
and
 the same exponent of SPT  for ABS if
 the problem $\widetilde{\B}$ is SPT for
ABS.
\end{rem}

\begin{rem}\label{rem2.2} Let the sequences $\az_j,\, \bz_j,\, \sz_j$ satisfy \eqref{2.002}.
If the problem $\B$  is  SPT (or PT) for ABS, then it is also SPT
(or PT) for NOR.

Indeed, it follows from \eqref{2.3} that $$\sum_{j=1}^\infty
\lz_{d,j}=\sum_{j=1}^d\az_j
+2\sum_{j=1}^d\bz_j\ze(\sz_j)\ge2\bz_1>0.$$ By \eqref{1.4} we have
$$n^{\B_d,\text{ABS}} (\va (2\bz_1)^{1/2})
 = \min \big\{ n\in \N : \sum_{j=n+1}^\infty\lz_{d,j}
\le2\bz_1\va^2 \le \va^2\sum_{j=1}^\infty \lz_{d,j}\big\},$$ which
means $$n^{\B_d,\text{NOR}}(\va )\le n^{\B_d,\text{ABS}}(\va
(2\bz_1)^{1/2}).$$  Therefore if the problem $\B$   is SPT (or PT)
for ABS, then it is also SPT (or PT) for NOR.

\end{rem}

\section{Proofs of Theorems 2.1-2.3}

\

\noindent{\it \textbf{Proof of Theorem \ref{th2.1}.}}

(i) From Remark \ref{rem2.2}, it is sufficient to show that $\B$
is PT for ABS.
 By Remark \ref{rem2.1}
we get that $\B$ is PT for ABS iff $\widetilde{\B}$ is PT for ABS.
Hence it suffices   to prove that $\widetilde{\B}$ is PT for ABS.

We note  $$\Big(\sum_{j=1}^\infty
\widetilde{\lz}_{d,j}^\tau\Big)^{\frac1{\tau}}
=\big(2\sum_{j=1}^d\bz_j^\tau\ze(\tau\sz_j)\big)^{\frac1{\tau}},$$
and \begin{equation}\label{2.8} 1\le \ze(\tau\sz_j)\le
\ze(\tau\sz_1)<+\infty,\quad\text{for}\quad\tau\in
(\frac1{\sz_1},1).
\end{equation}

It follows that for $\tau\in (\frac1{\sz_1},1)$,
\begin{align*}\Big(\sum_{j=1}^\infty \widetilde{\lz}_{d,j}^\tau\Big)^{\frac1{\tau}}
\le
\big(2\ze(\tau\sz_1)\big)^{\frac1{\tau}}\big(\sum_{j=1}^d\bz_j^\tau)^{\frac1{\tau}}
\le \big(2\ze(\tau\sz_1)\big)^{\frac1{\tau}}d^{\frac1{\tau}},
\end{align*}
where in the last inequality we used $\bz_j\le 1$ for all
$j\in\N$. This forces $$\sup_{d\in\N}\Big(\sum_{j=1}^\infty
\widetilde{\lz}_{d,j}^\tau\Big)^{\frac1{\tau}}d^{-\frac1{\tau}}
\le\big(2\ze(\tau\sz_1)\big)^{\frac1{\tau}}<+\infty.$$ Due to
\cite[Theorem 6.1]{NW1} we obtain that  $\widetilde{\B}$ is PT for
ABS. Therefore  the problem $\B$ is always PT for ABS or NOR.

(ii) From Remark \ref{rem2.1}, it is sufficient to prove that the problem
$$\widetilde{\B}=\{\widetilde{\B}_d:L_2([0,1]^d)\to L_2([0,1]^d)\},$$
 is SPT for ABS iff \eqref{2.4-1} holds, and that the exponent of
$\widetilde{\B}$ satisfies \eqref{2.4-2}.

Assume that \eqref{2.4-1} holds, i.e.,
$$A_*=\underset{d\to\infty}{\underline{\lim}}\frac{\ln\frac1{\bz_d}}{\ln
d}>1.$$ We want to prove that $\widetilde{\B}$ is SPT  for ABS.
Indeed, by \eqref{2.8} we have for any
$\tau\in(\frac1{\sz_1},1)$,
\begin{align}\label{3.17} \Big(\sum_{j=1}^\infty
\widetilde{\lz}_{d,j}^\tau\Big)^{\frac1{\tau}}&=
\Big(2\sum_{j=1}^d\bz_j^\tau\ze(\tau\sz_j)\Big)^{\frac1{\tau}}\notag\\
&\le\big(2\ze(\tau\sz_1)\big)^{\frac1{\tau}}\big(\sum_{j=1}^d\bz_j^\tau\big)^
{\frac1{\tau}}\notag\\ &\le
\big(2\ze(\tau\sz_1)\big)^{\frac1{\tau}}\big(\sum_{j=1}^\infty\bz_j^\tau\big)^
{\frac1{\tau}}.
\end{align}
Next, we shall prove that for any  $\tau\in(\frac1{A_*},1)$,
$$\big(\sum_{j=1}^\infty\bz_j^\tau\big)^
{\frac1{\tau}}<+\infty.$$ Indeed, for such
$\tau\in(\frac1{A_*},1)$, there exists a $\dz\in(0,1)$ for which
$$\tau=\frac1{A_*\,(1-\dz)}.$$   Since
$$A_*=\underset{d\to\infty}{\underline{\lim}}\frac{\ln\frac1{\bz_d}}{\ln
d}>1,$$  there exists a $d_0\in\N$ such that
 $$\frac{\ln\frac1{\bz_d}}{\ln d}\ge(1-\dz/2)\,A_*\quad\text{for}\quad d>d_0,$$
 which means $$\bz_d\le d^{-(1-\dz/2)A_*}\ \ {\rm for}\ \
 d>d_0.$$
 It follows that \begin{align*}
\sum_{j=1}^\infty\bz_j^\tau
\le\sum_{j=1}^{d_0}\bz_j^\tau+\sum_{d_0+1}^\infty
j^{-\frac{1-\dz/2}{1-\dz}}\le d_0+\sum_{j=1}^\infty
j^{-\frac{1-\dz/2}{1-\dz}}<+\infty,
\end{align*}
for any $\tau>\frac1{A_*}$. Hence we obtain
\begin{align}\label{3.000}\sup_{d\in\N}\Big(\sum_{j=1}^\infty
\widetilde{\lz}_{d,j}^\tau\Big)^{\frac1{\tau}}<+\infty,\end{align}
for any $\tau\in \big(\max\{\frac1{\sz_1}, \frac1{A_*}\},1\big)$.
We note that $\tau\in\big(\max\{\frac1{\sz_1},\frac1{A_*}\},1)$ is
equivalent to
$$\frac{2\tau}{1-\tau}> \frac2{\sz_1-1}\quad\text{and}
\quad\frac{2\tau}{1-\tau}> \frac2{A_*-1},$$ due to the
monotonicity of the function
$$\varphi(x)=\frac{x}{1-x},\ \,x\in(0,1).$$ It follows from \cite[Theorem 6.1]{NW1} that if \eqref{2.4-1} holds, then $\widetilde{\B}$ is SPT for
ABS,  and the exponent of SPT satisfies
\begin{equation}\label{3.18-1}p^{{\rm str-avg}}\le  \inf
\Big\{\frac{2\tau}{1-\tau}\ \big|\ \tau \ {\rm satisfies}\
\eqref{3.000}\Big\}\le\max\big\{\frac2{A_*-1},\frac2{\sz_1-1}\big\}.\end{equation}

On the other hand, assume that $\widetilde{\B}$ is SPT for ABS.
Then there exist positive $C,\, C_1,$ and $\tau\in(0,1)$ such that
\begin{align}\label{3.19}
+\infty>C&:=\sup_{d\in \N}\Big(\sum_{j=\lceil C_1\rceil}^\infty \widetilde{\lz}_{d,j}^\tau\Big)^{\frac1{\tau}}=\sup_{d\in \N}\Big(\sum_{j=1}^\infty
\widetilde{\lz}_{d,j}^\tau-\sum_{j=1}^{\lceil C_1\rceil-1}\widetilde{\lz}_{d,j}^\tau\Big)^{\frac1{\tau}}\notag\\
&=\sup_{d\in \N}\Big(2\sum_{j=1}^d\bz_j^\tau\ze(\tau\sz_j)
-\sum_{j=1}^{\lceil C_1\rceil-1}\widetilde{\lz}_{d,j}^\tau\Big)^{\frac1{\tau}}\notag\\
&\ge \sup_{d\in \N}\Big(2\sum_{j=1}^d\bz_j^\tau
-(\lceil C_1\rceil-1)\bz_1^\tau\Big)^{\frac1{\tau}}\notag\\
&\ge \sup_{d\in \N}\Big(2d\bz_d^\tau
-\lceil C_1\rceil+1\Big)^{\frac1{\tau}},
\end{align}
where  we used  $$1\ge \bz_1\ge\cdots >0,\ \ \
\widetilde{\lz}_{d,j}\le\bz_1\le1,\ \ \ze(\tau\sz_j)\ge1$$ in the
above inequalities.  Obviously, by \eqref{3.19} we have
$$\tau>\frac1{\sz_1}.$$
It follows from \eqref{3.19} that $$d\bz_d^\tau\le C_2,$$
 where $C_2:=\frac{C^\tau+\lceil C_1\rceil-1}{2}>0$,
which yields that  $$\frac{\ln\frac1{\bz_d}}{\ln
d}\ge\frac1{\tau}-\frac{\ln C_2}{\tau\ln d},$$ Letting $d\to
\infty$ we get
$$A_*=\underset{d\to\infty}{\underline{\lim}}\frac{\ln\frac1{\bz_d}}{\ln
d} \ge\frac1{\tau}>1.$$ which means $$\tau\ge\frac1{A_*}.$$

Hence if $\widetilde{\B}$ is SPT  for ABS, then  we have
 $$A_*=\underset{d\to\infty}{\underline{\lim}}\frac{\ln\frac1{\bz_d}}{\ln d}
>1,$$
and the exponent of SPT satisfies   $$p^{{\rm str-avg}}= \inf
\Big\{\frac{2\tau}{1-\tau}\ \big|\ \tau \ {\rm satisfies}\
\eqref{3.19}\Big\}\ge
\max\big\{\frac2{A_*-1},\frac2{\sz_1-1}\big\}.$$

Therefore ${\B}$ is SPT  for ABS iff
\begin{equation*}
A_*:=\underset{d\to\infty}{\underline{\lim}}\frac{\ln
\frac1{\bz_d}}{\ln d}>1,
\end{equation*}
 and the exponent of SPT  is
\begin{equation*}
p^{{\rm str-avg}}=\max\big\{\frac2{A_*-1},\frac2{\sz_1-1}\big\}.
\end{equation*}

Theorem \ref{th2.1} is proved. $\hfill\Box$

\

\noindent{\it \textbf{Proof of Theorem \ref{th2.2}.}}

From Remark \ref{rem2.2}, we note that  if $\B$ is SPT for ABS,
then it is SPT for NOR. If $A_*>1$ holds, by Theorem \ref{th2.1}
we get that  $\B$ is SPT for ABS, and hence is SPT for NOR. It
suffices to prove that if  $\B$ is SPT for NOR, then $A_*>1.$

Assume that $\B$ is SPT  for NOR. Then by   \cite[Theorem
6.2]{NW1} there exists a $\tau\in(0,1)$ such that
$$C:=\sup_{d\in \N}\frac{\Big(\sum_{j=1}^\infty \lz_{d,j}^\tau\Big)^{\frac1{\tau}}}{\sum_{j=1}^\infty \lz_{d,j}}
=\sup_{d\in
\N}\frac{\Big(\big(\sum_{j=1}^d\az_j\big)^\tau+2\sum_{j=1}^d\bz_j^\tau\ze(\tau\sz_j)\Big)^{\frac1{\tau}}}
{\sum_{j=1}^d\az_j+2\sum_{j=1}^d\bz_j\ze(\sz_j)}<+\infty.$$ It
follows that $$\tau>\frac1{\sz_1}.$$ Using \eqref{2.8}  and
$0\le\az_j\le c\bz_j$ for all $j\in \N$,
 we get
\begin{align}\label{3.7-1}
+\infty>C&\ge\frac{\Big(\big(\sum_{j=1}^d\az_j\big)^\tau+2\sum_{j=1}^d\bz_j^\tau\ze(\tau\sz_j)\Big)^{\frac1{\tau}}}
{\sum_{j=1}^d\az_j+2\sum_{j=1}^d\bz_j\ze(\sz_j)}\notag\\ &\ge
\frac{\Big(\big(\sum_{j=1}^d\az_j\big)^\tau+2\sum_{j=1}^d\bz_j^\tau\ze(\tau\sz_j)\Big)^{\frac1{\tau}}}
{(c+2\ze(\sz_1))\sum_{j=1}^d\bz_j}\notag\\
&\ge
\frac{2^{\frac1{\tau}}\big(\sum_{j=1}^d\bz_j^\tau\big)^{\frac1{\tau}}}{(c+2\ze(\sz_1))\sum_{j=1}^d\bz_j}.
\end{align}

Assume that $\sum_{j=1}^\infty\bz_j=+\infty$.  Then for the above
$\tau \in (0,1)$ we have
$$\big(\sum_{j=1}^\infty\bz_j^\tau\big)^{\frac1{\tau}}\ge\big(\sum_{j=1}^\infty\bz_j\big)^{\frac1{\tau}}=+\infty.$$
By the Stolz theorem we get
\begin{align}\label{3.8}
\lim_{d\to \infty}\frac{\big(\sum_{j=1}^d\bz_j^\tau\big)^{\frac1{\tau}}}{\sum_{j=1}^d\bz_j}&=\lim_{d\to \infty}\frac{\big(\sum_{j=1}^d\bz_j^\tau\big)^{\frac1{\tau}}-\big(\sum_{j=1}^{d-1}\bz_j^\tau\big)^{\frac1{\tau}}}{\bz_d}\notag\\
&=\lim_{d\to\infty}\frac{S_{\tau,d-1}\Big(\big(1+(\frac{\bz_d}{S_{\tau,d-1}})^\tau\big)^{\frac1{\tau}}-1\Big)}{\bz_d},
\end{align}
where
$$S_{\tau,d}=\big(\sum_{j=1}^d\bz_j^\tau\big)^{\frac1{\tau}}.$$

Since
$$0\le\frac{\bz_d}{S_{\tau,d-1}}\le\frac{\bz_1}{S_{\tau,d-1}}\to
0,\ \ {\rm as}\ \ d\to \infty,$$ and $${(1+x)^\az-1}\sim \az x, \ \
{\rm as}\ x\to 0$$  for  any  $\az>0,$ we have
\begin{align}\label{3.9}\Big(1+\big(\frac{\bz_d}{S_{\tau,d-1}}\big)^\tau\Big)^{\frac1{\tau}}-1
\sim\frac1{\tau}\big(\frac{\bz_d}{S_{\tau,d-1}}\big)^\tau,\quad
\text{as} \quad d\to \infty,\end{align} where $f(x)\sim g(x)$ as
$g(x)\to 0$ means that
$\lim\limits_{g(x)\to0}\frac{f(x)}{g(x)}=1$. Substituting
\eqref{3.9} into \eqref{3.8} yields
$$\lim_{d\to \infty}\frac{\big(\sum_{j=1}^d\bz_j^\tau\big)^{\frac1{\tau}}}{\sum_{j=1}^d\bz_j}
=\lim_{d\to\infty}\frac{S_{\tau,d-1}\cdot\frac1{\tau}\big(\frac{\bz_d}{S_{\tau,d-1}}\big)^\tau}{\bz_d}
=\lim_{d\to\infty}\frac1{\tau}\big(\frac{S_{\tau,d-1}}{\bz_d}\big)^{1-\tau}=+\infty,$$
contrary to \eqref{3.7-1}.

Hence we have $\sum_{j=1}^\infty\bz_j<+\infty$, and thence
$C_1:=\big(\sum_{j=1}^\infty\bz_j^\tau\big)^{\frac1{\tau}}<+\infty.$
It follows that
$$d\bz_d^\tau\le\sum_{j=1}^d\bz_j^\tau\le\sum_{j=1}^\infty\bz_j^\tau=C_1^\tau.$$
We obtain further  $$\frac{\ln\frac1{\bz_d}}{\ln
d}\ge\frac1{\tau}-\frac{\ln C_1}{\ln d}.$$ This gives
$$A_*=\underset{d\to\infty}{\underline{\lim}}\frac{\ln\frac1{\bz_d}}{\ln
d}\ge\frac1{\tau}>1,$$ which implies $\tau\ge\frac1{A_*}$. Hence
we conclude that if $\B$ is SPT   for NOR, then
\begin{equation*}A_*=\underset{d\to\infty}{\underline{\lim}}
\frac{\ln\frac1{\bz_d}}{\ln d}>1,\end{equation*} and the exponent
of SPT for NOR satisfies
 $$p^{{\rm str-avg}}\ge\max\big\{\frac2{A_*-1},\frac2{\sz_1-1}\big\}.$$

Theorem \ref{th2.2} is proved.
$\hfill\Box$

\

\noindent{\it \textbf{Proof of Theorem \ref{th2.3}.}}

First we assume that SPT holds for NOR.  By  \cite[Theorem
6.2]{NW1} there exists a $\tau\in(0,1)$ such that
$$C:=\sup_{d\in \N}\frac{\Big(\sum_{j=1}^\infty
\lz_{d,j}^\tau\Big)^{\frac1{\tau}}}{\sum_{j=1}^\infty
\lz_{d,j}}=\sup_{d\in
\N}\frac{\Big(\big(\sum_{j=1}^d\az_j\big)^\tau+2\sum_{j=1}^d\bz_j^\tau\ze(\tau\sz_j)\Big)^{\frac1{\tau}}}{\sum_{j=1}^d\az_j+2\sum_{j=1}^d\bz_j\ze(\sz_j)}<+\infty.$$
It follows that
$$\tau>\frac1{\sz_1}.$$
 We get for all $d\in\Bbb N$,
\begin{align*}+\infty>C&\ge\frac{\Big(\big(\sum_{j=1}^d\az_j\big)^\tau+2\sum_{j=1}^d\bz_j^\tau\ze(\tau\sz_j)\Big)^{\frac1{\tau}}}
{\sum_{j=1}^d\az_j+2\sum_{j=1}^d\bz_j\ze(\sz_j)}\ge \frac{2^{\frac1{\tau}}\big(\sum_{j=1}^d\bz_j^\tau\big)^{\frac1{\tau}}}{\sum_{j=1}^d\az_j+2\ze(\sz_1)\sum_{j=1}^d\bz_j}\\
&\ge\frac{\big(\sum_{j=1}^d\bz_j^\tau\big)^{\frac1{\tau}}}{\big(1+\frac{2\ze(\sz_1)}{r_1}\big)
\sum_{j=1}^d\az_j}
\ge\frac{d^{1/\tau}\bz_d}{c\big(1+\frac{2\ze(\sz_1)}{r_1}\big)d\az_d},
\end{align*}
where we used $$1\le \ze(t\sz_j)\le \ze(t\sz_1),\; \text{for
any}\; t\in(0,1],\;j\in\N,$$ and $1\ge\bz_1\ge\bz_2\ge\cdots\ge0$
in the last inequality. It follows that
$$\frac{\az_d}{\beta_d}\ge C_1d^{1/\tau -1}\quad\text{with}\quad
C_1:=\frac1{Cc\big(1+\frac{2\ze(\sz_1)}{r_1}\big)}.$$ This  means
$$\frac{\ln\frac{\az_d}{\bz_d}}{\ln d}\ge \frac1\tau -1+\frac{\ln C_1}{\ln d}.$$
Letting $d\to\infty$ we get
$$B_*=\underset{d\to\infty}{\underline{\lim}}\frac{\ln
\frac{\az_d}{\bz_d}}{\ln d}\ge\frac{1-\tau}{\tau}>0.$$ Hence if
$\B$ is SPT for NOR, then we have
$B_*=\underset{d\to\infty}{\underline{\lim}}\frac{\ln
\frac{\az_d}{\bz_d}}{\ln d}>0$, and the exponent of SPT for NOR
satisfies $$p^{{\rm
str-avg}}\ge\max\big\{\frac2{B_*},\frac2{\sz_1-1}\big\}.$$

On the other hand, we assume
$$B_*=\underset{d\to\infty}{\underline{\lim}}\frac{\ln\frac{\az_d}{\bz_d}}{\ln
d}>0.$$ Then for any
$\tau\in(\max\{\frac1{\sz_1},\frac1{1+B_*}\},1)$ we shall prove
$$\sup_{d\in \N}\frac{\Big(\sum_{j=1}^\infty
\lz_{d,j}^\tau\Big)^{\frac1{\tau}}}{\sum_{j=1}^\infty
\lz_{d,j}}<+\infty.$$ Indeed,  we have
\begin{align*}
\frac{\Big(\sum_{j=1}^\infty
\lz_{d,j}^\tau\Big)^{\frac1{\tau}}}{\sum_{j=1}^\infty \lz_{d,j}}&=
\frac{\Big(\big(\sum_{j=1}^d\az_j\big)^\tau+2\sum_{j=1}^d\bz_j^\tau\ze(\tau\sz_j)\Big)^{\frac1{\tau}}}
{\sum_{j=1}^d\az_j+2\sum_{j=1}^d\bz_j\ze(\sz_j)}\\
&\le\frac{2^{\frac1{\tau}}\Big(\sum_{j=1}^d\az_j+\big(2\ze(\tau\sz_1)\big)^{\frac1{\tau}}\big(\sum_{j=1}^d\bz_j^\tau\big)^
{\frac1{\tau}}\Big)}{\sum_{j=1}^d\az_j+2\sum_{j=1}^d\bz_j}\\
&\le
2^{\frac1{\tau}}+2^{\frac2{\tau}}\ze^{\frac1{\tau}}(\tau\sz_1)
\frac{\big(\sum_{j=1}^d\bz_j^\tau\big)^{\frac1{\tau}}}{\sum_{j=1}^d\az_j+2\sum_{j=1}^d\bz_j}
\\ &\le2^{\frac1{\tau}}+2^{\frac2{\tau}}\ze^{\frac1{\tau}}(\tau\sz_1)\frac{\big(\sum_{j=1}^d\bz_j^\tau\big)^{\frac1{\tau}}}{\sum_{j=1}^d\az_j},
\end{align*}
where in the first inequality we used
\begin{equation}\label{3.17-1} a+b\le \big(a^\tau+b^\tau
\big)^{\frac1{\tau}}\le 2^{\frac1{\tau}}(a+b),\quad a,\,b\ge0, \
\tau\in(0,1),\end{equation} and \eqref{2.8}. It suffices to prove
that
$$\sup_{d\in\N}\frac{\big(\sum_{j=1}^d\bz_j^\tau\big)^{\frac1{\tau}}}{\sum_{j=1}^d\az_j}<+\infty,\quad \text{for}\quad\tau\in(\frac1{1+B_*},1).$$

We suppose that $\tau=\frac1{1+B_*}\cdot\frac1{1-\dz_0}$ for some
$\dz_0\in(0,1)$. Since
$$B_*=\underset{d\to\infty}{\underline{\lim}}\frac{\ln\frac{\az_d}{\bz_d}}{\ln
d}=\underset{d\to\infty}{\underline{\lim}}\frac{\ln r_d}{\ln
d}>0,\ \ \  {\rm where}\ \ r_d=\frac{\az_d}{\bz_d},$$
 there exists a $d_0\in\N$ such that
$$\frac{\ln r_d}{\ln d}\ge(1-\dz_0)B_*=:\gz,\quad d>
d_0.$$It follows that  $$r_d^{-1}\le d^{-\gz}, \ \ d> d_0.$$ Note
that
$$\kappa:=\frac{\tau\gz}{1-\tau}=\frac{B_*-B_*\dz_0}{B_*-B_*\dz_0-\dz_0}>1.$$
Hence we get
\begin{align*}
\sum_{j=1}^\infty r_j^{-\frac{\tau}{1-\tau}}
&\le\sum_{j=1}^{d_0}r_j^{-\frac{\tau}{1-\tau}}+\sum_{j=d_0+1}^\infty j^{-\frac{\tau\gz}{1-\tau}}\\
&\le d_0r_1^{-\frac{\tau}{1-\tau}}+\sum_{j=1}^\infty
j^{-\kappa}<\infty .
\end{align*}
It follows from the H\"older inequality that
\begin{align*}&\quad \sup_{d\in\Bbb N}\frac{\big(\sum_{j=1}^d\bz_j^\tau\big)^{\frac1{\tau}}}{\sum_{j=1}^d\az_j}=
 \sup_{d\in\Bbb N} \frac{\big(\sum_{j=1}^d\az_j^\tau r_j^{-\tau}\big)^{\frac1{\tau}}}{\sum_{j=1}^d\az_j}
\\ &\le \sup_{d\in\Bbb N}\Big(\sum_{j=1}^d r_j^{-\frac \tau
{1-\tau}}\Big)^{\frac {1-\tau}\tau}  =\Big(\sum_{j=1}^\infty
r_j^{-\frac \tau {1-\tau}}\Big)^{\frac
{1-\tau}\tau}<\infty.\end{align*} Hence for any
$\tau\in(\max\{\frac1{\sz_1},\frac1{1+B_*}\},1)$  we have
$$\sup_{d\in \N}\frac{\Big(\sum_{j=1}^\infty \lz_{d,j}^\tau\Big)^{\frac1{\tau}}}{\sum_{j=1}^\infty \lz_{d,j}}<+\infty,$$
which means that SPT holds with the exponent of SPT for NOR
satisfying $$p^{{\rm
str-avg}}\le\max\big\{\frac2{B_*},\frac2{\sz_1-1}\big\}.$$

Therefore  SPT holds for NOR iff
$$B_*=\underset{d\to\infty}{\underline{\lim}}\frac{\ln\frac{\az_d}{\bz_d}}{\ln
d}>0$$ and  the exponent of SPT for NOR is $$p^{{\rm
str-avg}}=\max\big\{\frac2{B_*},\frac2{\sz_1-1}\big\}.$$

The proof of Theorem \ref{th2.2} is finished.
$\hfill\Box$

\

\section*{Acknowledgments}
 The  authors were supported
  by the National Natural Science
Foundation of China (Project no.  11671271),
 the  Beijing Natural Science Foundation (1172004).


\begin{thebibliography}{99}
\bibitem{CL} X. Chen, W. V. Li, Small deviation estimates for some additive processes,  Proceedings of the Conference on High Dimensional Probability III, Progress in Probability, 55 (2003) 225-238.
\bibitem{HWW} F. J. Hickernell, G. W. Wasilkowski, H. Wo\'zniakowski, Tractability of linear multivariate problems in the average-case
setting, in: A. Keller, S. Heinrich, H. Niederreiter (Eds.), Monte Carlo and Quasi-Monte Carlo Methods 2006, Springer,
Berlin, 2008, pp. 461-493.
\bibitem{KZ1} A. A. Khartov, M. Zani, Asymptotic analysis of average case approximation complexity of additive random fields, J. Complexity 52 (2019) 24-44.
\bibitem{KZ2} A. A. Khartov, M. Zani, Approximation complexity of sums of random processes, J. Complexity, in press.
\bibitem{LZ} M. A. Lifshitsa, M. Zani, Approximation complexity of additive random fields, J. Complexity  24 (2008) 362-379.
\bibitem{LZ2} M. A. Lifshits, M. Zani, Approximation of additive random
fields based on standard information: Average case and
probabilistic settings, J. Complexity 31 (5) (2015) 659-674.
\bibitem{NW1} E. Novak, H. Wo\'zniakowski, Tractability  of Multivariate Problems, Volume I: Liner Information, EMS, Z\"urich, 2008.
\bibitem{NW2} E. Novak, H. Wo\'zniakowski, Tractability  of Multivariate Problems, Volume II: Standard Information for Functionals, EMS, Z\"urich, 2010.
\bibitem{NW3} E. Novak, H. Wo\'zniakowski, Tractability  of Multivariate Problems, Volume III: Standard Information for Operators, EMS, Z\"urich, 2012.
\bibitem{S1} P. Siedlecki, Uniform weak tractability, J. Complexity 29 (6) (2013) 438-453.
\bibitem{WW} G. W. Wasilkowski, H. Wo\'zniakowski, Polynomial-time
algorithms for multivariate linear problems with finite-order
weights: average case setting, Found. Comput. Math. 9 (2009)
105-132.
\end{thebibliography}
\end{document}